\documentstyle[epsf,12pt]{article}
\def\Bbb R{{\rm \bf R}}
\def\proclaim#1{\vskip2mm{\bf #1}\em}
\def\endproclaim{\em \vskip2mm}
\def\tag#1{\eqno(#1)}
\def\gathered{\begin{array}{c}}
\def\endgathered{\end{array}}
\def\text{\mbox}

\begin{document}

\title {Probe method and a Carleman function}
\author{Masaru IKEHATA\\
Department of Mathematics, Graduate School of Engineering\\
Gunma University, Kiryu 376-8515, JAPAN}
\date{20 August 2007 Final}
\maketitle
\begin{abstract}
A Carleman function is a special fundamental solution with a
large parameter for the Laplace operator and gives a formula to
calculate the value of the solution of the Cauchy problem in a
domain for the Laplace equation. The probe method applied to an
inverse boundary value problem for the Laplace equation in a
bounded domain is based on the existence of a special sequence of
harmonic functions which is called a {\it needle sequence}. The
needle sequence blows up on a special curve which connects a given point
inside the domain with a point on the boundary of the domain and
is convergent locally outside the curve. The sequence yields a
reconstruction formula of unknown discontinuity, such as cavity,
inclusion in a given medium from the Dirichlet-to-Neumann map.
In this paper, an explicit needle sequence in {\it three dimensions}
is given in a closed form. It is an application of a Carleman
function introduced by Yarmukhamedov. Furthermore, an explicit
needle sequence in the probe method applied to the reduction of
inverse obstacle scattering problems with an {\it arbitrary} fixed wave
number to inverse boundary value problems for the Helmholtz equation is
also given.

\noindent
AMS: 35R30, 35R25, 35J05, 33E12, 35C05


\end{abstract}

\section{Introduction}
In \cite{Is} Isakov considered the uniqueness
issue of the problem for recovering the unknown inclusions in a
given electric conductive medium from infinitely many boundary
data.  This is a special, however, very important version of the Calder\'on
problem \cite{C} and is formulated as an inverse boundary value problem for an elliptic
partial differential equation in
which the boundary data are given by the associated Dirichlet-to-Neumann map.
Therein he established a uniqueness theorem.
The points of his paper are: a contradiction argument;
orthogonality identities
deduced by denying the conclusion, singular solutions, unique
continuation theorem, Runge approximation property and uniqueness of
the solution of the Cauchy problem for elliptic equations. See
\cite{IP} for these points.

After his work for almost 10 years nobody mentioned the {\it
reconstruction formula} of unknown inclusions. However,
finding such a formula is a natural attempt in studying inverse problems.
Nowadays we can cite at least three types formulae for the {\it full} reconstruction of unknown inclusions:

$\bullet$  The probe method \cite{Ik1} in two and three dimensions.

$\bullet$  An application by Br\"uhl \cite{Br} of the factorization method of Kirsch \cite{K1} in two and three dimensions.

$\bullet$  A solution of the Calder\'on problem in {\it two} dimensions by Astala-P\"aiv\"arinta \cite{AP1,AP2}.

This paper is concerned with the foundation of the probe method.
The method applied to inclusions in a medium with {\it constant}
conductivity is based on the existence of a special sequence of
harmonic functions which is called a {\it needle sequence} and
plays a role of a {\it probe needle}. The sequence blows up on a
special curve which is called a {\it needle} and connects a given
point inside the domain with a point on the boundary of the domain
and is convergent locally outside the curve. The {\it existence}
has been {\it ensured} by using the Runge approximation property
which can be proved by showing the denseness of the ranges of
infinitely many linear integral operators.  Therefore, the
construction of the needle sequence is not explicit and this makes
us difficult to understand theoretically the {\it point wise
behaviour} of the needle sequence on the needle.

In this paper, we give an explicit and concise needle sequence in
three dimensions for a special needle given by a line segment.
The idea for the construction came from reconsidering the role of a
Carleman function which is a special fundamental solution with a
large parameter for the Laplace operator and gives a formula to
calculate the value of the solution of the Cauchy problem in a
domain for the Laplace equation. Roughly speaking,  we say that a
function $\Phi(y,x,\tau)$ depending on a large parameter $\tau>0$
is called a Carleman function for a domain $\Omega$ and the
portion $\Gamma$ of $\partial\Omega$ if it satisfies the equation
$\triangle u+\delta(y-x)=0$ in $\Omega$; for each fixed
$x\in\Omega$ $\Phi(\,\cdot\,,x,\tau)$ and
$\partial/\partial\nu\Phi(\,\cdot\,,x,\tau)$ on
$\partial\Omega\setminus\Gamma$ vanish as
$\tau\longrightarrow\infty$. This function gives a representation
formula of any solution of the Cauchy problem for the Laplace
equation in $\Omega$ by using only the Cauchy data on $\Gamma$ and
yields a natural regularization for the numerical computation of
the solution from a noisy inaccurate Cauchy data on $\Gamma$. The
Cauchy problem for the Laplace equation is a fundamental and
important ill-posed problem appearing in mathematical sciences,
engineering and medicine. Therefore, it is quite important to seek
an explicit Carleman function in several domains. For this problem
Yarmukhamedov \cite{YM} gave a very interesting Carleman function
for a special domain in three dimensions which is a special
version of his fundamental solutions \cite{Y2} for the Laplace
operator.

In this paper, using his Carleman function, we give an {\it
explicit} needle sequence in three dimensions.  It
is written in a closed-form by using an integral involving Mittag-Leffler's function
$$\displaystyle
E_{\alpha}(z)=\sum_{n=0}^{\infty}\frac{z^n}{\Gamma(1+\alpha\,n)},\,\,0<\alpha\le 1.
$$
We never solve
integral equations in any sense.

\noindent We can summarize the conclusion of this paper in one
sentence:

$\bullet$  The {\it regular part} of a {\it special}
Carleman function which has been introduced
by Yarmukhamedov yields the desired needle sequence.
This makes an explicit {\it link} between the probe method and a Carleman function.

It should be pointed out that succeeding to the probe method, the author introduced another
method which he calls the {\it enclosure method} \cite{Ik2}. This
method gave us a different way of using the exponential solutions
for elliptic equations from Calder\'on 's way \cite{C} and yielded
an explicit extraction formula of the {\it convex hull} of the
unknown inclusions. This suggests that if one
replaces the exponential solutions with any other solutions having
a similar property, then one can get several information about
unknown inclusions.  The paper \cite{Ik9} is just a result in
which solutions coming from Mittag-Leffler's function
instead of the exponential solutions were employed.  The result
yielded more than the convex hull and the numerical implementation
has been done in \cite{ISS}.  However, this result is restricted
to the two dimensional case.  From this point of view the result
obtained in this paper can be considered as an extension of
\cite{Ik9} to three dimensions.

A brief outline of this paper is as follows.  In Section 2, we
explain what is the probe method and the role of the needle
sequence. For the purpose a simple inverse boundary value problem
for the Laplace equation is considered. Needless to say, a reader
who is familiar with the probe method can skip this section.  The
main idea is described in Section 3 by considering a two
dimensional case. In Section 4 we introduce Yarmukhamedov's
fundamental solution for the Laplace operator in three dimensions
and show that the regular part of a special version of his
fundamental solution yields an explicit needle sequence.  In
Section 5 we consider how to construct the needle sequence for the
Helmholtz equation.  We show that the Vekua transform of the
regular part of the special fundamental solution for the Laplace
operator yields the desired needle sequence. In Appendix, for
reader's convenience, we give a direct proof of a theorem
established by Yarmukhamedov.

\section{Probe method and needle sequence}

This section is devoted to a reader who is not familiar with the
probe method. Recently in \cite{Ik12, IkC, IkIm} the author
reformulated the probe method and further investigated the method
itself. It became much simpler than the previous formulation in
\cite{Ik1}.  In this section, for simplicity of description we
restrict ourselves to the case when the unknown discontinuity is
coming from {\it cavities}.

Let $\Omega$ be a bounded domain of $\Bbb R^m (m=2, 3)$ with a
Lipschitz boundary. Let $D$ be a bounded open set of $\Omega$ such
that $\overline D$ is contained in $\Omega$, $\partial D$ is
Lipschitz, $\Omega\setminus\overline D$ is connected.  Given $f\in
H^{1/2}(\partial\Omega)$ let $u\in H^1(\Omega\setminus\overline
D)$ be the unique weak solution of the problem
$$\begin{array}{c}
\displaystyle
\triangle u=0\,\,\text{in}\,\Omega\setminus\overline D,\\
\\
\displaystyle
\frac{\partial u}{\partial\nu}=0\,\,\text{on}\,\partial D,\\
\\
\displaystyle
u=f\,\,\text{on}\,\partial\Omega.
\end{array}
$$
The map $\Lambda_D:f\longmapsto \partial u/\partial\nu\vert_{\partial\Omega}$ is called the
Dirichlet-to-Neumann map.  We set $\Lambda_D=\Lambda_0$ in the
case when $D$ is empty.

In short, the probe method is a method of {\it probing inside}
given material by {\it monitoring} the behaviour of the sequence
of the energy gap
$$
\displaystyle
\int_{\partial\Omega}\{(\Lambda_0-\Lambda_{D})(v_n\vert_{\partial\Omega})\}\overline
v_n\vert_{\partial\Omega}\,dS
$$
for a specially chosen sequence $\{v_n\}$ of solutions of the
governing equation for the background medium ($D=\emptyset$) which
play a role of {\it probe needle}.

The method starts with introducing

{\bf\noindent Definition 2.1.}
Given a point $x\in\Omega$ we say that a non self-intersecting piecewise linear
curve $\sigma$ in $\overline\Omega$ is a needle with tip at $x$
if $\sigma$ connects a point on $\partial\Omega$
with $x$ and other points of $\sigma$ are contained in $\Omega$.

\noindent
We denote by $N_x$ the set of all needles with tip at $x$.

Let $\mbox{\boldmath $b$}$ be a nonzero vector in $\Bbb R^m$.
Given $x\in\Bbb R^m$, $\rho>0$ and $\theta\in]0,\pi[$ set
$$\displaystyle
C_x(\mbox{\boldmath $b$},\theta/2)
=\{y\in\,\Bbb R^m\,\vert\,
(y-x)\cdot\mbox{\boldmath $b$}>\vert y-x\vert\vert\mbox{\boldmath
$b$}\vert \cos(\theta/2)\}
$$
and
$$\displaystyle
B_{\rho}(x)=\{y\in\,\Bbb R^m\,\vert\,\vert y-x\vert<\rho\}.
$$
A set having the form
$$\displaystyle
V=B_{\rho}(x)\cap C_x(\mbox{\boldmath $b$},\theta/2)
$$
for some $\rho$, $\mbox{\boldmath $b$}$, $\theta$ and $x$
is called a {\it finite cone} with {\it vertex} at $x$.

Let $G(y)$ be a solution of the Laplace equation in
$\Bbb R^m\setminus\{0\}$ such that for any finite cone $V$ with vertex
at $0$
$$\displaystyle
\int_{V}\vert\nabla G(y)\vert^2 dy=\infty. \tag {2.1}
$$
Hereafter we fix this $G$.

For the new formulation of the probe method we need the following.

{\bf\noindent Definition 2.2.}
Let $\sigma\in N_x$.  We call the sequence $\xi=\{v_n\}$
of $H^1(\Omega)$ solutions of the Laplace equation
a {\it needle sequence} for $(x,\sigma)$ if it satisfies
for any compact set $K$ of $\Bbb R^m$ with $K\subset\Omega\setminus
\sigma$
$$\displaystyle
\lim_{n\longrightarrow\infty}(\Vert v_n(\,\cdot\,)-G(\,\cdot\,-x)\Vert_{L^2(K)}
+\Vert\nabla\{v_n(\,\cdot\,)-G(\,\cdot\,-x)\}\Vert_{L^2(K)})=0.
$$

The existence of the needle sequence is a consequence of
the Runge approximation property for the Laplace equation.

In \cite{Ik12} we clarified the behaviour of the needle sequence
on the needle as $n\longrightarrow\infty$.

\noindent
The two lemmas given below
are the core of the new formulation
of the probe method

\proclaim{\noindent Lemma 2.1.}
Let $x\in\Omega$ be an arbitrary point and $\sigma$ be a needle with tip at $x$.
Let $\xi=\{v_n\}$ be an arbitrary needle sequence for $(x,\sigma)$.
Then, for any finite cone $V$ with vertex at $x$ we have
$$\displaystyle
\lim_{n\longrightarrow\infty}\int_{V\cap\Omega}\vert\nabla v_n(y)\vert^2dy=\infty.
$$
\endproclaim

\proclaim{\noindent Lemma 2.2.}
Let $x\in\Omega$ be an arbitrary point and $\sigma$ be a needle with tip at $x$.
Let $\xi=\{v_n\}$ be an arbitrary needle sequence for $(x,\sigma)$.
Then for any point $z\in\sigma$ and open ball $B$ centered at $z$ we have
$$\displaystyle
\lim_{n\longrightarrow\infty}\int_{B\cap\Omega}\vert\nabla v_n(y)\vert^2dy=\infty.
$$
\endproclaim

\noindent
From these lemmas we know that one can recover {\it full knowledge} of the given needle
as the set of all points where the needle sequence for the needle blows up.
This means that the needle is realized as a special sequence of harmonic functions
without loosing information about the geometry of the needle.  This is the new point
added on the probe method in \cite{Ik12}.

\noindent
{\bf\noindent Definition 2.3.}
Given $x\in\,\Omega$, needle $\sigma$ with tip $x$
and needle sequence $\xi=\{v_n\}$ for $(x,\sigma)$
define
$$\displaystyle
I(x,\sigma,\xi)_n=\int_{\partial\Omega}\{(\Lambda_0-\Lambda_{D})\overline f_n\}
f_n\,dS,\,n=1,2,\cdots
$$
where
$$\displaystyle
f_n(y)=v_n(y),\,\,y\in\,\partial\Omega.
$$

$\{I(x,\sigma,\xi)_n\}_{n=1,2,\cdots}$ is a sequence
depending on $\xi$ and $\sigma\in N_x$.
We call the sequence the {\it indicator sequence}.

The behaviour of the indicator sequence has two sides.
One side is closely related to the function defined below.

\noindent
{\bf\noindent Definition 2.4.}
The {\it indicator function} $I$ is defined by the formula
$$\displaystyle
I(x)=\int_D\vert\nabla G(y-x)\vert^2 dy+\int_{\Omega\setminus\overline D}
\vert\nabla w_x\vert^2 dy,\,\,x\in\Omega\setminus\overline D
$$
where $w_x$ is the unique weak solution of the problem:
$$\begin{array}{c}
\displaystyle
\triangle w=0\,\,\text{in}\,\Omega\setminus\overline D,\\
\\
\displaystyle
\frac{\partial w}{\partial\nu}=-\frac{\partial}{\partial\nu}(G(\,\cdot\,-x))\,\,\text{on}\,\partial D,\\
\\
\displaystyle w=0\,\,\text{on}\,\partial\Omega.
\end{array}
$$

The function $w_x$ is called the {\it reflected solution} by $D$.

The following theorem says that

$\bullet$ one can calculate the value of the indicator function at an arbitrary point
outside the cavity from $\Lambda_0-\Lambda_D$;

$\bullet$ the indicator function can not be continued across the boundary of the cavity
as a bounded function in the whole domain.

\proclaim{\noindent Theorem 2.1.A}
We have

$\bullet$  (A.1)  given $x\in\Omega\setminus\overline D$ and needle $\sigma$ with tip at $x$
if $\sigma\cap\overline D=\emptyset$,
then for any needle sequence $\xi=\{v_n\}$ for $(x,\sigma)$
the sequence $\{I(x,\sigma,\xi)_n\}$ converges to the indicator function $I(x)$

$\bullet$  (A.2)  for each $\epsilon>0$
$$\displaystyle
\sup_{\displaystyle\text{dist}\,(x,\,D)>\epsilon}I(x)<\infty
$$

$\bullet$  (A.3)  for any point $a\in\,\partial D$
$$\displaystyle
\lim_{x\longrightarrow a}I(x)=\infty.
$$

\endproclaim

Since mathematically Theorem 2.1.A is enough for establishing a
reconstruction formula of the cavities, in the previous applications of the probe
method \cite{Ik1} we did not consider the natural question

$\bullet$ what happens on the indicator sequence when the tip of the
needle is just located on the boundary of cavities, inside or passing
through the cavities?

\noindent
However, in practice the tip of the needle cannot move forward with
infinitely small step
and therefore in the scanning process with needle there is a
possibility of skipping the unknown boundary of cavities, entering
inside or passing through the cavities. So for the
practical use of the probe method we have to clarify the behaviour
of the indicator sequence in those cases.  The answer to this
question is

\proclaim{\noindent Theorem 2.1.B}
Let $x\in\Omega$ and $\sigma\in N_x$.
If $x\in\Omega\setminus\overline D$ and
$\sigma\cap D\not=\emptyset$
or $x\in\overline D$, then for any needle sequence $\xi=\{v_n\}$ for $(x,\sigma)$ we have
$\lim_{n\longrightarrow\infty}I(x,\sigma,\xi)_n=\infty$.
\endproclaim

\begin{figure}[htbp]
\begin{center}
\epsfxsize=10cm
\epsfysize=10cm
\epsfbox{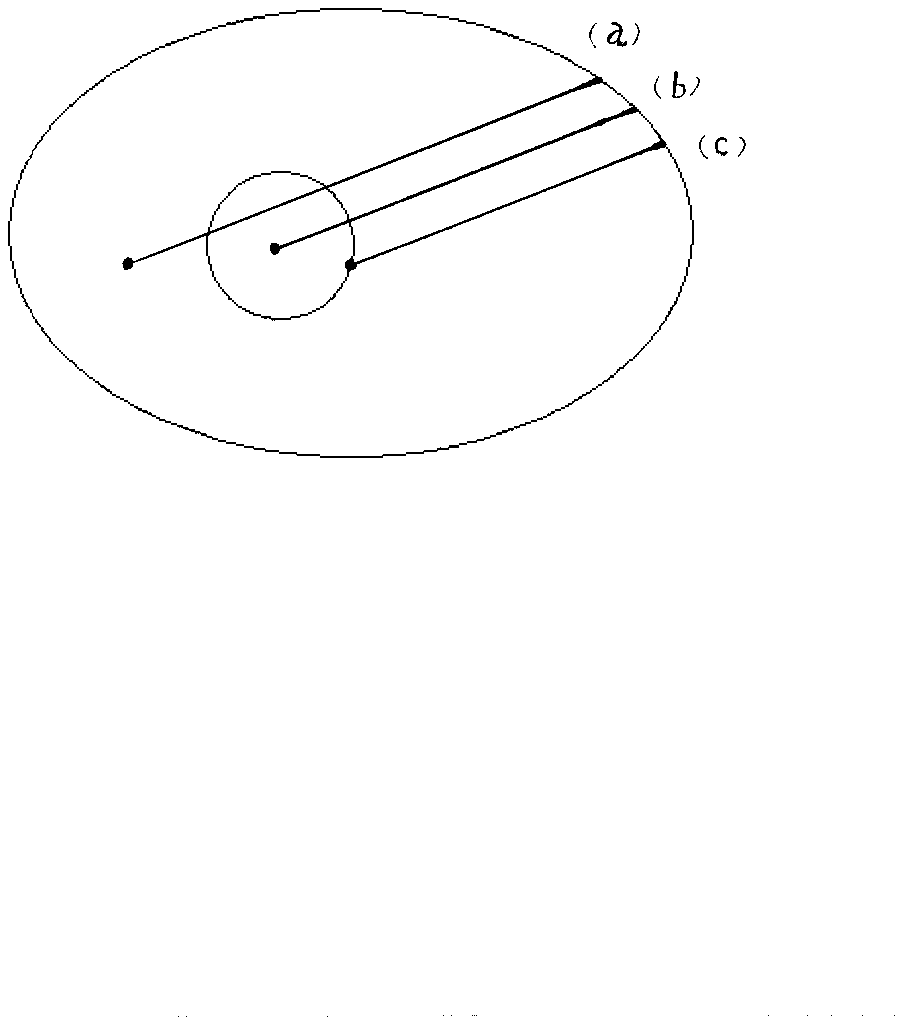}
\caption{
(a) $x\in\Omega\setminus\overline D$ and $\sigma\cap D\not=\emptyset$.
(b) $x\in D$.  (c) $x\in\partial D$.
}\label{fig1}
\end{center}
\end{figure}

\noindent
See Figure \ref{fig1} for typical situations.

\vspace{-1.0cm}

These two theorems are essentially a special case of the results in \cite{Ik12}.
Finally we note that, as a corollary of
Theorems 2.1 A and B one gets a characterization of $\overline D$.

\proclaim{\noindent Corollary 2.4}
A point $x\in\,\Omega$ belongs to $\Omega\setminus\overline D$ if and only if
there exists a needle $\sigma$ with tip at $x$ and needle sequence $\xi$
for $(x,\sigma)$ such that the indicator sequence $\{I(x,\sigma,\xi)_n\}$ is bounded
from above.

\endproclaim

\noindent
Needless to say, this automatically gives a uniqueness theorem too.

\section{Main idea:  a new role of a Carleman function}

\noindent We think that the reader now understand the basic of the probe method and
the role of the needle sequence.
So the next problem is {\it how to construct} the needle sequence.

There are several points for the meaning of the `construction'. It
is well known that one can construct the needle sequence by
solving {\it infinitely} many first kind integral equations in the
sense of {\it minimum norm solutions}. The minimum norm solutions
are given by a combination of the Tikhonov {\it regularization method}
and Morozov {\it discrepancy principle}(see \cite{KRS} for these concepts).
However, from our point of
view this does not yield an explicit needle sequence.  One can
also point out that the determination of the regularization
parameter via the Morozov discrepancy principle itself is a
nonlinear problem and seems impossible to find the parameter
explicitly.

So how can one construct the needle sequence {\it explicitly}?
In this and following sections we always consider
a geometrically simplest needle described in

{\bf\noindent Definition 3.1}
A needle with tip at $x$ is called a {\it straight needle} with tip at $x$ directed to $\omega$
if the needle is given by $l_x(\omega)\cap\overline\Omega$ where
$$
\displaystyle
l_{x}(\omega)=\{x+t\omega\,\vert\,0\le t<\infty\}.
$$

\noindent In two dimensions in \cite{IN} we constructed an explicit needle sequence
for $G$ given by
$$\displaystyle
G(y)=\frac{1}{y_1+iy_2}.
\tag {3.1}
$$
Note that this $G$ satisfies (2.1) for any finite cone $V$ with vertex
at $0$.

For the construction we made use of a well known idea in complex function theory:
$$\begin{array}{c}
\displaystyle
\frac{1}{y-x}
=\frac{1}{(y-x')-(x-x')}\\
\\
\displaystyle
=\frac{1}{\displaystyle (y-x')\left(1-\frac{x-x'}{y-x'}\right)}
\\
\\
\displaystyle
=\sum_{m=0}^{\infty}\frac{(x-x')^m}{(y-x')^{m+1}}
\end{array}
$$
provided $x'\approx x$ and $x'$ is located on a needle with tip at $x$.  Note that we used the same symbol
for the complex number $z_1+iz_2$ corresponding to a point $z=(z_1,z_2)^T$.

However, in general we have to do this type of procedure many times and
therefore the resulted sequence is quite involved.

In this section, we give a different and extremely simple needle
sequence and explain the idea behind.

{\bf\noindent Definition 3.2.}
Given a unit vector $\omega=(\omega_1,\,\omega_2)^T$ in $\Bbb R^2$ and $\alpha\in\,]0,\,1]$
define
$$\displaystyle
v(y;\alpha,\tau,\omega)=-\frac{E_{\alpha}(\tau\,y\,\overline\omega)-1}{y},\,\,\tau>0.
$$
This is a harmonic function in the whole plane.
Note also that $\displaystyle y\cdot(\omega+i\omega^{\perp})=y\,\overline\omega$
where $\omega^{\perp}=(-\omega_2,\omega_1)$.

\proclaim{\noindent Theorem 3.1}
Let $x\in\Omega$ and $\sigma$ be
a straight needle with tip at $x$ directed to $\omega$.
The sequence $\{v(\,\cdot\,-x;\alpha_n,\tau_n,\omega)\vert_{\Omega}\}$ is a needle sequence for
$(x,\,\sigma)$ where $\alpha_n$ and
$\tau_n$ are suitably chosen sequences and satisfy

$\bullet$  $0<\alpha_n\le 1$, $\tau_n>0$

$\bullet$  $\alpha_n\longrightarrow 0$ and $\tau_n\longrightarrow\infty$ as $n\longrightarrow\infty$.

\endproclaim

{\it\noindent Proof.}
For each fixed $\alpha$ we have, as $\tau\longrightarrow\infty$,
$$\displaystyle
v(\,\cdot\,-x;\alpha,\tau,\omega)\longrightarrow G(\,\cdot\,-x)
$$
uniformly for $y$ in any compact subset of $\Bbb R^2\setminus\overline{C_x(\pi\alpha/2,\omega)}$.
This is a direct consequence of the well known asymptotic behaviour
of Mittag-Leffler's function \cite{B, E}: as $\vert z\vert\longrightarrow\infty$
with $\pi\alpha/2+\epsilon\le\vert\text{arg}\,z\vert\le\pi$
$$\displaystyle
E_{\alpha}(z)=-\frac{1}{z\,\Gamma(1-\alpha)}+O\left(\frac{1}{\vert z\vert^2}\right)
\tag {3.2}
$$
where $\alpha$ and $\epsilon$ are fixed and satisfies $0<\alpha<1$ and $0<\epsilon<\pi-\pi\alpha/2$.

The choice of $\{\alpha_n\}$ and $\{\tau_n\}$ depend on that of an
{\it exhaustion} for $\Omega\setminus\sigma$.  More
precisely let $\{O_n\}$ be a sequence of open subsets of
$\Omega\setminus\sigma$ such that for all $\overline
O_n\subset O_{n+1}$; $\cup_n O_n=\Omega\setminus\sigma$.
We call this sequence an exhaustion for
$\Omega\setminus\sigma$.   Let $\{\epsilon_n\}$ be an
arbitrary sequence of positive numbers such that
$\epsilon_n\longrightarrow 0$ as $n\longrightarrow\infty$. The
choice of $\{\alpha_n\}$ and $\{\tau_n\}$ can be done as follows.
Since $\overline O_1\subset\Omega\setminus\sigma$, one
can choose a small $\alpha\in\,]0,\,1[$ in such a way that
$\overline O_1$ is contained in the set $\Bbb
R^2\setminus\overline {C_x(\pi\alpha/2,\omega)}$. Set
$\alpha_1=\alpha$. Since $v(y-x;\alpha_1,\tau,\omega)$ converges to
$G(y-x)$ in $H^1(\Omega_1)$ as $\tau\longrightarrow\infty$, one
can find a large $\tau$ in such a way that $\Vert
v(\,\cdot\,-x;\alpha_1,\tau,\omega)-G(\,\cdot\,-x)\Vert_{H^1(O_1)}<\epsilon_1$.
Set $\tau_1=\tau$. For $\alpha_2,\cdots$ and $\tau_2,\cdots$ we do
a similar procedure.  Then the obtained sequence satisfies the
desired property.

\noindent
$\Box$

\noindent
We see that $v(y-x;\alpha_n,\tau_n,\omega)$ takes the special value at the tip $x$ of the needle:
$$\displaystyle
v(y-x;\alpha_n,\tau_n,\omega)\vert_{y=x}
=-\frac{\tau_n\,\overline\omega}
{\Gamma(1+\alpha_n)}
\tag {3.3}
$$
and thus this yields
$$\displaystyle
\lim_{n\longrightarrow\infty}\frac{v(y-x;\alpha_n,\tau_n,\omega)\vert_{y=x}}{\tau_n}=-\overline\omega.
$$
Therefore the leading term of the asymptotic behaviour of this
special needle sequence at the tip $x$ of needle $\sigma$ does not
depend on the choice of $\{\alpha_n\}$ and is uniquely determined
by $\{\tau_n\}$.

However, the point of this section is not the statement of
Theorem 3.1 itself.  The main point is how to find this simple form.
For simplicity let $x=0$.

First the function $G(y)$ given by (3.1) is a fundamental solution for the operator $\overline\partial$
(ignoring multiplying a constant, hereafter same).  Since $E_{\alpha}(0)=1$, the function
$$\displaystyle
\frac{E_{\alpha}(\tau\,y\,\overline\omega)}{y}
\tag {3.4}
$$
is also a fundamental solution for operator $\overline\partial$.
The point is: this fundamental solution becomes small outside a sector as
$\tau\longrightarrow\infty$ \cite{B, E}.  Since we can write
$$\displaystyle
\frac{E_{\alpha}(\tau y\,\overline\omega)}{y}
=G(y)+\frac{E_{\alpha}(\tau y\,\overline\omega))-1}{y},
$$
this means that the outside the sector
$$\displaystyle
G(y)\approx -\frac{E_{\alpha}(\tau y\,\overline\omega)-1}{y}.
$$
Note that this right-hand side is just coming from the regular part
of (3.4).

The function (3.4) has been appeared in Yarmukhamedov's work
\cite{Y6} to give an explicit formula of the value of the solution
of the Cauchy problem in a special domain for $\overline\partial$.  Such a type of
fundamental solution for some operator is called the Carleman
function (for the operator) and the formula is called a Carleman-type formula.
See also \cite{A} for Carleman-type formulae in complex function theory
and \cite{AB, IC} for the stationary Schr\"odinger equation.

So now the principle of the choice of needle sequence is clear.
In particular, in three dimensions
the principle is: {\it choose $(-1)$ times the {\it regular part} of a suitable
Carleman function (for the Laplace operator)}.

\noindent This principle is the conclusion of this paper and gives
a {\it new role} for the Carleman function. We found a direct link
between the Carleman function and the probe method.
Needless to
say, it is different from the idea in existing author's previous
application \cite{Iu} of the Carleman function to an inverse
boundary value problem in unbounded domain.

In the next section we consider the three dimensional case more precisely.

\section{Yarmukhamedov's fundamental solution and needle sequence in three dimensions}

The following is taken from \cite{Y2}; however, therein
the proof is not given.  For reader's convenience we gave a direct proof of this
theorem in Appendix.

\proclaim{\noindent Theorem 4.1}
Let $\lambda\ge 0$.
Let $K(w)$ be an entire function such that

$\bullet$  $K(w)$ is real for real $w$

$\bullet$  $K(0)=1$

$\bullet$  for each $R>0$  and $m=0, 1, 2$
$$\displaystyle
\sup_{\vert\text{Re}\,w\vert<R}\vert K^{(m)}(w)\vert e^{\lambda\vert\text{Im}\,w\vert}
<\infty.
\tag {4.1}
$$
Define
$$\displaystyle
-2\pi^2\Phi(x)
=\frac{1}{2}\int_0^{\infty}\text{Im}\,(\frac{K(w)}{w})\frac{e^{\lambda u}+e^{-\lambda u}}
{\sqrt{\vert x'\vert^2+u^2}}du.
\tag {4.2}
$$
where $w=x_3+i\sqrt{\vert x'\vert^2+u^2}$ and $x'=(x_1,x_2)$.
Then one has the expression
$$\displaystyle
\Phi(x)=\frac{e^{i\lambda\vert x\vert}}{4\pi\vert x\vert}+H(x)
$$
where $H$ is $C^2$ in the whole space and satisfies
$$\displaystyle
\triangle H(x)+\lambda^2 H(x)=0\,\,\text{in}\,\Bbb R^3
$$
and therefore $\Phi$ satisfies
$$\displaystyle
\triangle\Phi(x)+\lambda^2\Phi(x)+\delta(x)=0\,\,\text{in}\,\Bbb R^3.
$$
\endproclaim

Hereafter we set $\Phi=\Phi_K$ to denote the dependence
on $K$ and consider only the case when $\lambda=0$.
In this case (3.2) ensures that $K(w)=E_{\alpha}(\tau w)$ with $\tau>0$ satisfies (4.1).
Yarmukhamedov \cite{YM} established that the function $\Phi_K$
for this $K$ with a fixed $\alpha$ is a Carleman function (for the Laplace equation).

{\bf\noindent Definition 4.1}
Given two unit vectors $\vartheta_1$ and $\vartheta_2$ in three dimensions
and $\alpha\,\in]0,\,1]$ define
$$\displaystyle
v(y;\alpha,\tau,\vartheta_1,\vartheta_2)
=-\left\{\Phi_{K}(y\cdot\vartheta_1,y\cdot\vartheta_2,y\cdot(\vartheta_1\times\vartheta_2))
-\frac{1}{4\pi\,\vert y\vert}\right\},\,\,\tau>0
$$
where $K(w)=E_{\alpha}(\tau w)$.  From Theorem 4.1 for $\lambda=0$ one knows that this function of $y$ is
nothing but the regular part
of $\Phi_K$ and harmonic in the whole space.

In this section we prove

\proclaim{\noindent Theorem 4.2} Let $x\in\Omega$ and $\sigma$ be
a straight needle with tip at $x$ directed to
$\omega=\vartheta_1\times\vartheta_2$. Then the sequence
$\{v(\,\cdot\,-x;\alpha_n,\tau_n,\vartheta_1,\vartheta_2)\vert_{\Omega}\}$ is a
needle sequence for $(x,\sigma)$ with
$$\displaystyle
G(y)=\frac{1}{4\pi\vert y\vert}
$$
where $\alpha_n$ and $\tau_n$ are suitably chosen sequences and satisfy

$\bullet$  $0<\alpha_n<1$, $\tau_n>0$

$\bullet$ $\alpha_n\longrightarrow 0$ and $\tau_n\longrightarrow\infty$.

\endproclaim

{\it\noindent Proof.}
First we prove that, for each fixed $\alpha$ as $\tau\longrightarrow\infty$
$$\displaystyle
v(\,\cdot\,-x;\alpha,\tau,\vartheta_1,\vartheta_2)\longrightarrow G(\,\cdot\,-x)
$$
uniformly for $y$ in any compact subset of $\Bbb R^3\setminus\overline{C_x(\pi\alpha/2,\omega)}$.

Given $\eta\in]0, \pi[$, $r>0$ let $\gamma(\eta,r)$ denote the
contour that originates at $\infty$, runs in toward the origin
just above the half line $z=\rho e^{-i\eta}, \rho\ge r$, the part
of the circle $z=re^{i\theta}, \vert\theta\vert\le\eta$, the origin
counterclockwise and then returns to $\infty$ just above the half
line $z=\rho e^{i\eta}, \rho\ge r$. The contour $\gamma$ splits
the complex plane into the two simply connected infinite domains
$D^{-}(\eta,r)$ and $D^{+}(\eta,r)$ lying, respectively to the left and the right
of $\gamma(\eta,r)$.

From Mittag-Leffler's integral representation (p. 206 of \cite{B}) and analytic continuation,
it follows that, for any $\eta\in]\pi/2, \pi[$
$$\displaystyle
E_{\alpha}(z)=\frac{1}{2\pi i}
\int_{\gamma(\eta,r)}\frac{\zeta^{\alpha-1}e^{\zeta}}{\zeta^{\alpha}-z}d\zeta,\,\,z\in\,D^{-}(\alpha\eta,r^{\alpha}).
\tag {4.3}
$$

\noindent
Let $\tau>0$.  We see that: if $z\in\,D^{-}(\alpha\eta,
r^{\alpha})\setminus B_{r^{\alpha}}(0)$, then $\tau z\in
D^{-}(\alpha\eta, r^{\alpha})$; if $z\in\,D^{-}(\alpha\eta,
r^{\alpha})$, then so is $\overline z$. And it should be pointed
out that $E_{\alpha}(z)$ is real for real $z$. From these we have,
for all $z\in\,D^{-}(\alpha\eta, r^{\alpha})\setminus
B_{r^{\alpha}}(0)$,
$$\displaystyle
\frac{1}{\text{Im}\,z}\text{Im}\,\left(\frac{E_{\alpha}(\tau z)}{z}\right)
=-\frac{1}{2\pi\,i\vert z\vert^2}
\int_{\gamma(\eta,r)}
\frac{(\zeta^{\alpha}-2\tau\,\text{Re}\,z)\zeta^{\alpha-1}e^{\zeta}}{(\zeta^{\alpha}-\tau\,z)
(\zeta^{\alpha}-\tau\,\overline z)}d\zeta.
\tag {4.4}
$$

Here we restrict the interval where $\eta$ belongs to
$\eta\in\,]\pi/2,\,\pi/(2\alpha)[$. Let $\theta$ satisfy
$\alpha\eta<\theta<\pi/2$.  Then we have, for all $w\in
D^{-}(\theta, r^{\alpha})\setminus B_{r^{\alpha}}(0)$ and $\tau>1$,
$$\displaystyle
\inf_{\zeta\in\gamma(\eta,r)}\vert \zeta^{\alpha}-\tau w\vert
\ge C \tau\vert w\vert
\tag {4.5}
$$
where
$$
C=\min\,\{1-\frac{1}{\tau}, \sin\,(\theta-\alpha\eta)\}.
$$
This can be proved as follows.  It suffices to consider the case when
$\theta\le \text{arg}\,w\equiv\varphi\le\pi$.  We divide the case into two subcases.
First let $\theta\le \varphi\le\alpha\eta+\pi/2$.  In this case
the circle $\vert z-w\vert=\vert w\vert\sin(\varphi-\alpha\eta)$ has a single common point
with the half line $z=ke^{i\alpha\eta},\,r^{\alpha}\le k<\infty$.  Thus this yields
$$\displaystyle
\inf_{\zeta\in\gamma(\eta,r)}\vert \zeta^{\alpha}-w\vert=\vert w\vert \sin\,(\varphi-\alpha\eta).
$$
Since $\theta-\alpha\eta\le\varphi-\alpha\eta\le\pi/2$, we obtain
$$
\displaystyle
\inf_{\zeta\in\gamma(\eta,r)}\vert \zeta^{\alpha}-w\vert\ge \vert w\vert \sin\,(\theta-\alpha\eta).
\tag {4.6}
$$

Next consider the case $\alpha\eta+\pi/2<\varphi\le\pi$.  A simple geometrical observation yields
$$\displaystyle
\inf_{\zeta\in\gamma(\eta,r)}\vert \zeta^{\alpha}-w\vert
=\vert w-r^{\alpha}e^{i\alpha\eta}\vert
\ge\vert w\vert-r^{\alpha}.
\tag {4.7}
$$
Now let $\tau>1$.  Then $\tau w$ still belongs to $D^{-}(\theta,r^{\alpha})\setminus B_{r^{\alpha}}(0)$
and having the same arg as $w$.  Therefore, inserting $\tau w$ instead of $w$ into (4.7) we obtain
$$\displaystyle
\inf_{\zeta\in\gamma(\eta,r)}\vert \zeta^{\alpha}-\tau w\vert
\ge \tau\vert w\vert-r^{\alpha}
\ge (\tau-1)\vert w\vert.
\tag {4.8}
$$

\noindent
Now a combination of (4.6) and (4.8) gives (4.5).

\noindent
Now let $y$ satisfy
$$\displaystyle
(y-x)\cdot\omega\le\vert y-x\vert\cos\theta
\tag {4.9}
$$
and
$$\displaystyle
r^{\alpha}\le\vert y-x\vert\le R^{\alpha}.
\tag {4.10}
$$
Set $w=(y-x)\cdot\omega+i\sqrt{\vert (y-x)\cdot\vartheta_1\vert^2+
\vert (y-x)\cdot\vartheta_2\vert^2+u^2}$.  This $w$ belongs to
$D^{-}(\theta,r^{\alpha})\setminus B_{r^{\alpha}}(0)$.  Therefore from (4.4) and (4.5) we obtain
$$\displaystyle
\vert \frac{1}{\text{Im}\,w}\text{Im}\,\left(\frac{E_{\alpha}(\tau w)}{w}\right)\vert
\le \frac{C}{2\pi\,\tau\,(r^{2\alpha}+u^2)^2}
\int_{\gamma(\eta,r)}(\vert\zeta\vert^{\alpha}+2R^{\alpha})\vert\zeta\vert^{\alpha-1}
e^{\displaystyle\text{Re}\,\zeta}\vert d\zeta\vert.
$$
The integral in this right hand side is absolutely convergent
since $\pi/2<\eta<\pi$. Therefore, we have
$$\displaystyle
\Phi_{K}((y-x)\cdot\vartheta_1, (y-x)\cdot\vartheta_2,
(y-x)\cdot\omega)=O(\frac{1}{\tau})
$$
as $\tau\longrightarrow\infty$ and the estimate is uniform for
$y-x$ with (4.9) and (4.10) for each fixed $\theta$, $\eta$, $r$,
$R$ and $\alpha$.

It is easy to see that, choosing suitable $\theta$, $\eta$, $r$, a
large fixed $R$ and $\alpha$, one can generate an exhaustion of
the domain $\Omega\setminus\sigma$. Therefore from the
principle explained in Section 3 we conclude that
the sequence
$\{v(\,\cdot\,-x;\alpha_n,\tau_n,\vartheta_1,\vartheta_2)\vert_{\Omega}\}$ converges to
$G(y)$ in $L^2_{\text{loc}}(\Omega\setminus\sigma)$ for suitably chosen
sequences $\{\alpha_n\}$ and $\{\tau_n\}$.  Then, a standard argument gives
the convergence in $H^1_{\text{loc}}(\Omega\setminus\sigma)$.

\noindent
$\Box$

{\bf\noindent Remark 4.1.}
The proof presented here can be considered as a minor modification of that of Lemma in \cite{YM}.
Anyway the new point that should be emphasized
is the {\it relationship} between the construction of a needle sequence
and a Carleman function.

{\bf\noindent Definition 4.2.}
We call the needle sequence $\{v(\,\cdot\,-x;\alpha_n,\tau_n,\vartheta_1,\vartheta_2)\vert_{\Omega}\}$ given
in Theorem 4.2 a {\it standard} needle sequence for $(x,\sigma)$ for the Laplace equation.

Choosing $K(w)=1$ in $\Phi_K$, we have
$$\displaystyle
\Phi_1(y)=\frac{1}{4\pi\vert y\vert}.
$$
This gives the expression:
$$\displaystyle
2\pi^2 v(y;\alpha,\tau,\vartheta_1,\vartheta_2)
=\int_0^{\infty}\text{Im}\,\left(\frac{E_{\alpha}(\tau\,w)-1}{w}\right)\,\frac{du}
{\sqrt{\vert y\cdot\vartheta_1\vert^2
+\vert y\cdot\vartheta_2\vert^2+u^2}}
\tag {4.11}
$$
where
$$\displaystyle
w=y\cdot\omega+i
{\sqrt{\vert y\cdot\vartheta_1\vert^2
+\vert y\cdot\vartheta_2\vert^2+u^2}}.
$$

\noindent
The advantage of the standard needle sequence is that one can exactly
know the precise values in a closed-form on the needle and its linear extension.

\proclaim{\noindent Theorem 4.3} Let
$\omega=\vartheta_1\times\vartheta_2$.
We have:

(1) if $y=x+s\,\omega$ with $s\not=0$, then
$$\displaystyle
v(y-x;\alpha,\tau,\vartheta_1,\vartheta_2)
=\frac{E_{\alpha}(\tau\,s)-1}
{4\pi\,s};
\tag {4.12}
$$

(2) if $y=x$, then
$$\displaystyle
v(y-x;\alpha,\tau,\vartheta_1,\vartheta_2)\vert_{y=x}=\frac{\tau}{4\pi\Gamma(1+\alpha)}.
\tag {4.13}
$$

\endproclaim

{\it\noindent Proof.}

\noindent
From (4.11) we have the expression of $v(y-x;\alpha,\tau,\vartheta_1,\vartheta_2)$
in both cases (1) and (2):
$$\displaystyle
2\pi^2 v(y-x;\alpha,\tau,\vartheta_1,\vartheta_2)
=\int_0^{\infty}\text{Im}\,\left(\frac{E_{\alpha}(\tau(x_3+iu))-1}{x_3+iu}\right)\,\frac{du}{u}
\tag {4.14}
$$
here we set $x_3=(y-x)\cdot\omega$ for simplicity of description.
\noindent
One can write
$$\displaystyle
\int_0^{\infty}\text{Im}\,\left(\frac{E_{\alpha}(\tau(x_3+iu))-1}{x_3+iu}\right)\,\frac{du}
{u}
=\frac{1}{2i}\lim_{\epsilon\downarrow 0}
\left(\int_{-i/\epsilon}^{-i\epsilon}+\int_{i\epsilon}^{i/\epsilon}\right)
\frac{E_{\alpha}(\tau(x_3+\zeta))-1}{(x_3+\zeta)\zeta}\,d\zeta.
\tag {4.15}
$$

\noindent
First consider the case when $x_3\not=0$.  Then $\zeta=-x_3$ is a removable singular point
for the integrand.  Therefore, from the Cauchy's integral theorem and (3.2) we see that
the right hand side of (4.15) becomes
$$
\displaystyle
\frac{1}{2i}\lim_{\epsilon\downarrow 0}
\int_{\zeta=\epsilon e^{i\theta},\,
\vert\theta-\pi\vert\le\pi/2}
\frac{E_{\alpha}(\tau(x_3+\zeta))-1}{(x_3+\zeta)\zeta}\,d\zeta
=\frac{\pi}{2}\,
\frac{E_{\alpha}(\tau\,x_3)-1}{x_3}.
$$
This together with (4.14) and $(y-x)\cdot\omega=s$ for $y=x+s\,\omega$
yields (4.12).

A similar computation in the case when $x_3=0$ yields (4.13).

\noindent
$\Box$

The formula (4.13) corresponds to (3.3).

{\bf\noindent Remark 4.2.}
From (4.12), (4.13) and the trivial estimate
$$\displaystyle
E_{\alpha}(x)-1>\frac{x^n}{\Gamma(1+\alpha\,n)}, \,x>0,\,n=1,2,\cdots
$$
we see that the standard needle sequence on the needle blows up in
a {\it point wise} sense:  for all points $y$ on $l_{x}(\omega)$
(see Definition 3.1) $\lim_{n\longrightarrow\infty}
v(y-x;\alpha_n,\tau_n,\vartheta_1,\vartheta_2)=+\infty$.

{\bf\noindent Remark 4.3.}
We can know also the values of $\nabla v(y-x;\alpha,\tau,\vartheta_1,\vartheta_2)$ on
the line $y=x+s\,\omega\,(-\infty<s<\infty)$ with $\omega=\vartheta_1\times\vartheta_2$.
Define
$$\displaystyle
I_1(\rho,s)=\int_1^{\infty}
\text{Im}\,\left(\frac{E_{\alpha}(\tau(s+i\sqrt{\rho+u^2}))}{s+i\sqrt{\rho+u^2}}\right)
\frac{du}{\sqrt{\rho+u^2}},\,\rho>-1,\,-\infty<s<\infty
$$
and
$$
\displaystyle
I_2(\rho,s)
=\int_0^1
\text{Im}\,\left(\frac{E_{\alpha}(\tau(s+i\sqrt{\rho+u^2}))}{s+i\sqrt{\rho+u^2}}\right)
\frac{du}{\sqrt{\rho+u^2}},\,0<\rho<1,\,-\infty<s<\infty.
$$
Using (3.2) and the asymptotic expansions of the derivatives of
Mittag-Leffler's function outside the sector
$\vert\text{arg}\,z\vert\le\pi\alpha/2$, we conclude that
$I_1(\rho,s)$ is a smooth function of $(\rho,s)$.  Following the
argument in Appendix A.2 we have
$$\begin{array}{c}
\displaystyle
I_2(\rho,s)
=-\int_0^1\frac{du}{\rho+s^2+u^2}\\
\\
\displaystyle
+\sum_{n=1}^{\infty}
\frac{\tau^{n+1}}{\Gamma(1+\alpha(n+1))}\,
\sum_{\text{$j$ is odd}}
\frac{n!}{j!(n-j)!}s^{n-j}
(-1)^{(j-1)/2}
\int_0^1(\rho+u^2)^{(j-1)/2}du.
\end{array}
$$
Since the second term is smooth for $(\rho,s)\,\in\Bbb R^2$, this yields that the function
$$\displaystyle
I_3(\rho,s)\equiv I_2(\rho,s)+\int_0^1\frac{du}{\rho+s^2+u^2}
$$
has a smooth extension to $\Bbb R^2$.
Since
$$\displaystyle
\int_0^{\infty}\text{Im}\,\left(\frac{1}{s+i\sqrt{\rho+u^2}}\right)
\frac{du}{\sqrt{\rho+u^2}}
=-\int_0^{\infty}\frac{du}{\rho+s^2+u^2},
$$
we have the decomposition of (4.11):
$$
\displaystyle
2\pi^2v(y;\alpha,\tau,\vartheta_1,\vartheta_2)
=I_1(\rho,s)+I_3(\rho,s)
\tag {4.16}
$$
where $\rho=\vert y\cdot\vartheta_1\vert^2+\vert y\cdot\vartheta_2\vert^2$ and
$s=y\cdot\omega$.
Therefore we obtain
$$\displaystyle
2\pi^2\nabla v(y;\alpha,\tau,\vartheta_1,\vartheta_2)
=\frac{\partial}{\partial s}(I_1+I_3)\,\omega
+2\frac{\partial}{\partial\rho}(I_1+I_3)\,
\{(y\cdot\vartheta_1)\vartheta_1+(y\cdot\vartheta_2)\vartheta_2\}.
$$
This together with (4.12) and (4.16) yields that for $y=x+s\,\omega$ with $s\not=0$
$$\displaystyle
\nabla v(y-x;\alpha,\tau,\vartheta_1,\vartheta_2)
=\frac{d}{ds}\left\{\frac{E_{\alpha}(\tau s)-1}{4\pi\,s}\right\}\omega.
$$
Moreover, letting $s\longrightarrow 0$, we have
$$
\displaystyle \nabla
v(y-x;\alpha,\tau,\vartheta_1,\vartheta_2)\vert_{y=x}
=\frac{\tau^2}{4\pi\Gamma(1+2\alpha)}\,\omega.
$$
Therefore $\nabla v$ on the line
$y=x+s\,\omega\,(-\infty<s<\infty)$ is {\it parallel} to the
direction of the line.  In particular, using the power series
expansion of Mittag-Leffler's function, we see that $\nabla
v(y-x;\alpha,\tau,\vartheta_1,\vartheta_2)\cdot\omega>0$ and
$\lim_{n\longrightarrow\infty}\nabla
v(y-x;\alpha_n,\tau_n,\vartheta_1,\vartheta_2)\cdot\omega=+\infty$
on the set $l_x(\omega)$.

\section{Needle sequence for Helmholtz equation and Vekua transform}

\noindent
The needle sequence for the Helmoltz equation can be defined as same as
the Laplace equation.

Let $G(y)$ be a solution of the Helmholtz equation
$\triangle v+\lambda^2v=0$ in $\Bbb R^3\setminus\{0\}$ such that
for any finite cone $V$ with vertex at $0$ $G$ satisfies (2.1).
Hereafter we fix $G$.

{\bf\noindent Definition 5.1.}
Let $\sigma\in N_x$.  We call the sequence $\xi=\{v_n\}$
of $H^1(\Omega)$ solutions of the Helmholtz equation
$\triangle v+\lambda^2 v=0$
a {\it needle sequence} for $(x,\sigma)$ for the Helmholtz equation if it satisfies,
for any compact set $K$ of $\Bbb R^3$ with $K\subset\Omega\setminus
\sigma$
$$\displaystyle
\lim_{n\longrightarrow\infty}(\Vert v_n(\,\cdot\,)-G(\,\cdot\,-x)\Vert_{L^2(K)}
+\Vert\nabla\{v_n(\,\cdot\,)-G(\,\cdot\,-x)\}\Vert_{L^2(K)})=0.
$$

{\bf\noindent Remark 5.1.}
Any needle sequences for the Helmholtz equation with an arbitrary
$\lambda$ blows up on the needle in the sense of Lemmas 2.1 and 2.2.
See \cite{Ik12} for the proof.

\noindent
The existence of the needle sequence for the Helmholtz equation has been
ensured under the additional condition on $\lambda$
(see Appendixes of \cite{IW,Ik12} for the proof):

$\bullet$ $\lambda^2$ is not a Dirichlet eigenvalue
of $-\triangle$ in $\Omega$.

\noindent
However, the proof therein does not give us any explicit form of a needle sequence.
Since the needle sequence for the Helmholtz equation
plays a central role in the probe method applied to inverse obstacle scattering problems
(see \cite{IW,Ik12} again and \cite{IkIm} for a recent application),
it is quite important to obtain an explicit one.

In this section we give explicit needle sequences for the
Helmholtz equation for all straight needles and all $\lambda(>0)$.
There is no additional condition on $\lambda$.
For the construction we do not make use of Theorem 4.1.
Mittag-Leffler's function with $\alpha<1$ does not satisfy the
condition (4.1) when $\lambda>0$. Therefore one can not substitute
$K(w)=E_{\alpha}(\tau\,w)$ into (4.2) to construct a needle sequence for
the Helmoltz equation $\triangle u+\lambda^2 u=0$ in three
dimensions. In this section, instead of (4.2) we employ an idea of
making use of a transformation introduced by Vekua \cite{V, V2}.

Let $\lambda>0$. The Vekua
transform $v\longmapsto T_{\lambda}v$ in three dimensions takes
the form
$$\displaystyle
T_{\lambda}v(y)=v(y)-\frac{\lambda\vert y\vert}{2}\,\int_0^1v(ty)J_1(\lambda\vert y\vert\sqrt{1-t})
\,\sqrt{\frac{t}{1-t}}\,dt
$$
where $J_1$ stands for the Bessel function of order $1$.
The important property of this transform
is: if $v$ is harmonic in the whole space, then $T_{\lambda}v$ is a solution of the Helmholtz
equation $\triangle u+\lambda^2 u=0$ in the whole space.

{\bf\noindent Definition 5.2.}
Given two unit vectors $\vartheta_1$ and $\vartheta_2$ in three dimensions
and $\alpha\,\in]0,\,1]$
define
$$\displaystyle
v^{\lambda}(y;\alpha,\tau,\vartheta_1,\vartheta_2)
=T_{\lambda}v(y;\alpha,\tau,\vartheta_1,\vartheta_2),\,\,\tau>0.
$$
This function of $y$ satisfies the Helmholtz equation
$\triangle u+\lambda^2 u=0$ in the whole space.

\proclaim{\noindent Theorem 5.1}  Let $x\in\Omega$ and $\sigma$ be a straight needle with tip at $x$ directed to
$\omega=\vartheta_1\times\vartheta_2$.
Then the sequence $\{v^{\lambda}(\,\cdot\,-x;\alpha_n,\tau_n,\vartheta_1,\vartheta_2)\vert_{\Omega}\}$
is a needle sequence for $(x,\sigma)$ for the Helmholtz equation
with $G=G_{\lambda}$ given by
$$\displaystyle
G_{\lambda}(y)=\text{Re}\,\left(\frac{e^{i\lambda\vert y\vert}}{4\pi\vert y\vert}\right)
$$
where $\alpha_n$ and $\tau_n$ are suitably chosen sequences and satisfy

$\bullet$  $0<\alpha_n<1$, $\tau_n>0$

$\bullet$ $\alpha_n\longrightarrow 0$ and $\tau_n\longrightarrow\infty$.

\endproclaim

{\it\noindent Proof.}
It suffices to prove the theorem in the case when $x=0$,
$\vartheta_1=(1,0,0)^T$, $\vartheta_2=(0,1,0)^T$ and thus
$\omega=(0,0,1)^T$.

First we show that for each fixed $\alpha$ as
$\tau\longrightarrow\infty$ the function
$v^{\lambda}(y;\alpha,\tau,\vartheta_1,\vartheta_2)$ converges to
the function $v(y)$ defined by
$$\displaystyle
v(y)=\Phi_1(y)-\frac{\lambda\vert y\vert}{2}\,\int_0^1
\Phi_1(ty)J_1(\lambda\vert y\vert\sqrt{1-t})
\sqrt{\frac{t}{1-t}}dt,\, y\not=0
\tag {5.1}
$$
uniformly for $y$ in any compact subset of $\Bbb R^3\setminus\overline{C_0(\pi\alpha/2,\omega)}$.

From Theorem 4.2 we have, $\Phi_K(y)$ with $K(z)=E_{\alpha}(\tau\,z)$
converges to $0$ uniformly for $y$ in any compact subset of $\Bbb R^3\setminus\overline{C_{0}(\pi\alpha/2,\omega)}$
as $\tau\longrightarrow\infty$.
Since $\Phi_1(y)-v(y;\alpha,\tau,\vartheta_1,\vartheta_2)=\Phi_{K}(y)$,
it suffices to prove that the following function of $y$
$$\displaystyle
R(y;\alpha,\tau)=\int_0^1\Phi_{K}(ty)\,J_1(\lambda\vert
y\vert\sqrt{1-t})\,\sqrt{\frac{t}{1-t}}\,dt,\,\,y\not=0,
\tag {5.2}
$$
converges to $0$ uniformly for $y$ in any compact subset of $\Bbb R^3\setminus\overline{C_{0}(\pi\alpha/2,\omega)}$.

Let $\eta\in\,]\pi/2,\pi[$. Let $r$ and $R$ be arbitrary positive numbers with $r<R$.
Let $\theta\in\,]\alpha\eta,\,\pi/2[$.
Let $y$ satisfy (4.9) and (4.10) with $x=0$ and $\omega=(0,0,1)^T$.
Then the non zero complex number
$$\displaystyle
w=y_3+i\sqrt{y_1^2+y_2^2+u^2}, u\ge 0
$$
satisfies $\theta\le\vert\text{arg}\,w\vert\le\pi$. Then we see
that, for all $\tau>0$ $\tau\, w\in\, D(\alpha\eta,r^{\alpha})$.
Thus (4.3) gives
$$\displaystyle
E_{\alpha}(\tau w)
=\frac{1}{2\pi\,i}
\int_{\gamma(\eta,r)}
\frac{\zeta^{\alpha-1}\,e^{\zeta}}
{\zeta^{\alpha}-\tau\,w}\,d\zeta.
$$
Since $\overline{E_{\alpha}(z)}=E_{\alpha}(\overline z)$ for all $z$, we obtain
$$\begin{array}{c}
\displaystyle
-2\pi^2\Phi_K(y)
=\int_0^{\infty}\frac{1}{\text{Im}\,w}\text{Im}\,\left(\frac{E_{\alpha}(\tau\,w)}
{w}\right)du\\
\\
\displaystyle
=-\frac{1}{2\pi i}\int_0^{\infty}\left\{\frac{1}{\vert w\vert^2}
\int_{\gamma(\eta,r)}
\frac{(\zeta^{\alpha}-2\tau\text{Re}\,w)\zeta^{\alpha-1}e^{\zeta}}
{(\zeta^{\alpha}-\tau w)(\zeta^{\alpha}-\tau\overline w)}\,d\zeta\right\}du.
\end{array}
$$
Let $0<t<1$.  A change of variable $u\longrightarrow tu$ yields
$$\displaystyle
-2\pi^2\Phi_K(ty)
=-\frac{1}{2\pi i t}\int_0^{\infty}\left\{\frac{1}{\vert y\vert^2+u^2}
\int_{\gamma(\eta,r)}
\frac{(\zeta^{\alpha}-2\tau\,ty_3)\zeta^{\alpha-1}e^{\zeta}}
{(\zeta^{\alpha}-\tau\,t\, w)(\zeta^{\alpha}-\tau\,t\,\overline w)}\,d\zeta\right\}du.
\tag {5.3}
$$
Divide $\gamma(\eta,r)$ as $\gamma_1+\gamma_2+\gamma_3$ where
$$\begin{array}{c}
\displaystyle
\gamma_1: \zeta=se^{i\eta},\,r\le s<\infty,\\
\\
\displaystyle
-\gamma_2: \zeta=se^{-i\eta},\,\,r\le s<\infty,\\
\\
\displaystyle
\gamma_3: \zeta=re^{i\varphi},\,\,\vert\varphi\vert\le\eta.
\end{array}
$$
Then from (5.2) and (5.3) we have the expression
$$\displaystyle
4\pi^3 i R(y;\alpha_n,\tau_n)\equiv R_1+R_2+R_3
$$
where
$$
\displaystyle
R_j
=\int_0^1\left(\int_0^{\infty}\left\{\frac{1}{\vert y\vert^2+u^2}
\int_{\gamma_j}
\frac{(\zeta^{\alpha}-2\tau\,ty_3)\zeta^{\alpha-1}\,e^{\zeta}}
{(\zeta^{\alpha}-\tau\,t\, w)(\zeta^{\alpha}-\tau\,t\,\overline w)}\,d\zeta\right\}du
\frac{J_1(\lambda\vert y\vert\sqrt{1-t})}{\sqrt{t(1-t)}}\right)dt.
\tag {5.4}
$$

\noindent
From a simple geometrical observation,
for all $t\in\,[0,\,1]$ and $u\ge 0$ we have the estimates
$\displaystyle\vert e^{i\alpha\eta}-tw\vert\ge\sin\,(\theta-\alpha\eta)$
and $\displaystyle
\vert e^{i\alpha\eta}-t\overline w\vert\ge\sin\,(\theta-\alpha\eta)$.
Note that the common right-hand side of theses estimates gives the distance
between the point $e^{i\alpha\eta}$ and the set
$\{w\,\vert\,\theta\le\vert\text{arg}\,w\vert\le\pi\}$.

Then these estimates together with a change of variable yield
$$\begin{array}{c}
\displaystyle
\vert\int_{\gamma_1}
\frac{(\zeta^{\alpha}-2\tau\,ty_3)\zeta^{\alpha-1}\,e^{\zeta}}
{(\zeta^{\alpha}-\tau\,t\, w)(\zeta^{\alpha}-\tau\,t\,\overline w)}\,d\zeta\vert
\le\int_{r}^{\infty}
\frac{s^{\alpha-1}e^{s\cos\,\eta}(s^{\alpha}+2\tau\,\vert y_3\vert)ds}
{\vert s^{\alpha}e^{i\alpha\,\eta}-\tau\,t\,w\vert
\vert s^{\alpha}e^{i\alpha\,\eta}-\tau\,t\,\overline w\vert}\\
\\
\displaystyle
=\frac{1}{\alpha}\int_{r^{\alpha}/\tau}^{\infty}
\frac{e^{\tau^{1/\alpha}\xi^{1/\alpha}\cos\,\eta}(\xi+2\vert y_3\vert)d\xi}
{\vert e^{i\alpha\,\eta}-t\,w\vert
\vert e^{i\alpha\,\eta}-t\,\overline w\vert}\\
\\
\displaystyle
\le
\frac{1}{\alpha\,\sin^2(\theta-\alpha\eta)}\int_{0}^{\infty}
e^{\tau^{1/\alpha}\xi^{1/\alpha}\cos\,\eta}(\xi+2\vert y_3\vert)d\xi\\
\\
\displaystyle
\le
\frac{1}{\alpha\,\sin^2(\theta-\alpha\eta)\tau^{1/\alpha}}
\int_0^{\infty}e^{s\cos\,\eta}(\tau^{-1/\alpha}s+2R^{\alpha})ds
=O(\tau^{-1/\alpha}).
\end{array}
$$
Similarly we have the same estimate for the integral of the same
function on $\gamma_2$.  Therefore, using the well known estimate
$$
\vert J_1(s)\vert\le\frac{s}{2},\,s\ge 0,
$$
we obtain, for $j=1,2$ $R_j=O(\tau^{-1/\alpha})$ as $\tau\longrightarrow\infty$.
This is a uniform estimate with $y$ satisfying (4.9) and (4.10) for $x=0$ and $\omega=(0,0,1)^T$.

Next we show that $R_3=O(\tau^{-1/2})$ as $\tau\longrightarrow\infty$.

\noindent
We have
$$\begin{array}{c}
\displaystyle
I(w,\tau t)\equiv \vert\int_{\gamma_3}
\frac{(\zeta^{\alpha}-2\tau\,ty_3)\zeta^{\alpha-1}\,e^{\zeta}}
{(\zeta^{\alpha}-\tau\,t\, w)(\zeta^{\alpha}-\tau\,t\,\overline w)}\,d\zeta\vert\\
\\
\displaystyle
\le
r^{\alpha}
\int_{-\eta}^{\eta}e^{r\,\cos\,\varphi}
\frac{\vert r^{\alpha}e^{i\alpha\varphi}-2\tau t\,y_3\vert d\varphi}
{\vert r^{\alpha}e^{i\alpha\varphi}-\tau\,t\,w\vert
\vert r^{\alpha}e^{i\alpha\varphi}-\tau\,t\,\overline w\vert}\\
\\
\displaystyle
\le
(rR)^{\alpha}(1+2\tau\,t)
\int_{-\eta}^{\eta}
\frac{d\varphi}
{\vert r^{\alpha}e^{i\alpha\varphi}-\tau\,t\,w\vert
\vert r^{\alpha}e^{i\alpha\varphi}-\tau\,t\,\overline w\vert}
\end{array}
$$
and thus this yields
$$\begin{array}{c}
\displaystyle
\vert R_3\vert\le\int_0^1\left(\int_0^{\infty}\frac{I(w,\tau t)\,du}{\vert y\vert^2+u^2}\right)
\frac{\vert J_1(\lambda\vert y\vert\sqrt{1-t})\vert dt}{\sqrt{t(1-t)}}
\\
\\
\displaystyle
\le
(rR)^{\alpha}\frac{\lambda R^{\alpha}}{2}
\int_0^1\frac{1+2\tau\,t}{\sqrt{t}}\,dt
\int_0^{\infty}\frac{du}{\vert y\vert^2+u^2}
\left(\int_{-\eta}^{\eta}
\frac{d\varphi}
{\vert r^{\alpha}e^{i\alpha\varphi}-\tau\,t\,w\vert
\vert r^{\alpha}e^{i\alpha\varphi}-\tau\,t\,\overline w\vert}\right).
\end{array}
\tag {5.5}
$$

A change of variable gives
$$\begin{array}{c}
\displaystyle
\int_0^1\frac{1+2\tau\,t}{\sqrt{t}}\,dt
\int_0^{\infty}\frac{du}{\vert y\vert^2+u^2}
\left(\int_{-\eta}^{\eta}
\frac{d\varphi}
{\vert r^{\alpha}e^{i\alpha\varphi}-\tau\,t\,w\vert
\vert r^{\alpha}e^{i\alpha\varphi}-\tau\,t\,\overline w\vert}\right)
\\
\\
\displaystyle
=\frac{1}{\sqrt{\tau}}
\int_0^{\tau}\frac{1+2\xi}{\sqrt{\xi}}d\xi
\int_0^{\infty}\frac{du}{\vert y\vert^2+u^2}
\left(\int_{-\eta}^{\eta}
\frac{d\varphi}
{\vert r^{\alpha}e^{i\alpha\varphi}-\xi\,w\vert
\vert r^{\alpha}e^{i\alpha\varphi}-\xi\,\overline w\vert}\right)
\\
\\
\displaystyle
\le\frac{1}{\sqrt{\tau}}
\int_0^{\infty}\frac{1+2\xi}{\sqrt{\xi}}d\xi
\int_0^{\infty}\frac{du}{\vert y\vert^2+u^2}
\left(\int_{-\eta}^{\eta}
\frac{d\varphi}
{\vert r^{\alpha}e^{i\alpha\varphi}-\xi\,w\vert
\vert r^{\alpha}e^{i\alpha\varphi}-\xi\,\overline w\vert}\right)
=\frac{1}{\sqrt{\tau}}(I_1+I_2)
\end{array}
\tag {5.6}
$$
where
$$\displaystyle
I_1=\int_2^{\infty}\frac{1+2\xi}{\sqrt{\xi}}d\xi
\int_0^{\infty}\frac{du}{\vert y\vert^2+u^2}
\left(\int_{-\eta}^{\eta}
\frac{d\varphi}
{\vert r^{\alpha}e^{i\alpha\varphi}-\xi\,w\vert
\vert r^{\alpha}e^{i\alpha\varphi}-\xi\,\overline w\vert}\right)
$$
and
$$\displaystyle
I_2=\int_0^{2}\frac{1+2\xi}{\sqrt{\xi}}d\xi
\int_0^{\infty}\frac{du}{\vert y\vert^2+u^2}
\left(\int_{-\eta}^{\eta}
\frac{d\varphi}
{\vert r^{\alpha}e^{i\alpha\varphi}-\xi\,w\vert
\vert r^{\alpha}e^{i\alpha\varphi}-\xi\,\overline w\vert}\right).
$$

\noindent
If $\xi\ge 2$, then we have
$$
\displaystyle
\vert r^{\alpha}e^{i\alpha\varphi}-\xi\,w\vert
\ge\xi\vert w\vert-r^{\alpha}
\ge\xi\vert w\vert(1-\frac{r^{\alpha}}{\xi\vert w\vert})
\ge\xi\vert w\vert (1-\frac{r^{\alpha}}{2\vert y\vert})
\ge\frac{\xi\vert w\vert}{2}.
$$
Therefore we get
$$\displaystyle
I_1\le
4\int_2^{\infty}\frac{1+2\xi}{\xi^{2+1/2}}d\xi\int_0^{\infty}\frac{du}{(r^{2\alpha}+u^2)^2}
<\infty.
\tag {5.7}
$$

On the other hand, a geometrical observation gives, for all
$\xi>0$
$$\displaystyle
\vert r^{\alpha}e^{i\alpha\varphi}-\xi w\vert
\ge r^{\alpha}\sin\,(\theta-\alpha\eta).
$$
This right-hand side is just the distance between two sets
$\{w\,\vert\, \theta\le\vert\text{arg}\,w\vert\le\pi\}$
and $\gamma_3$ in the complex plane.
Therefore we obtain
$$\displaystyle
I_2\le\frac{2\eta}{r^{2\alpha}\sin^2\,(\theta-\alpha\eta)}
\int_0^2\frac{1+2\xi}{\sqrt{\xi}}\,d\xi
\int_0^{\infty}\frac{du}{r^{2\alpha}+u^2} <\infty.
\tag {5.8}
$$

\noindent
Now from (5.5) to (5.8) we obtain the desired estimate for $R_3$.

Summing up, we have $R=O(\tau^{-1/2})$ as $\tau\longrightarrow\infty$
and this is a uniform estimate with $y$ satisfying (4.9) and (4.10) for $x=0$ and $\omega=(0,0,1)^T$.
Therefore the function $v^{\lambda}(\,\cdot\,;\alpha,\tau,\vartheta_1,\vartheta_2)$
converges to the function $v(y)$ uniformly for $y$ in any
compact subset of $\Bbb R^3\setminus\overline{C_{0}(\pi\alpha/2,\omega)}$.
Since $\alpha$ can be arbitrary small and the function
$v^{\lambda}(\,\cdot\,;\alpha,\tau,\vartheta_1,\vartheta_2)$
satisfies the elliptic equation $\triangle u+\lambda^2
u=0$, a standard argument yields that,
for suitably chosen sequences $\{\alpha_n\}$ and $\{\tau_n\}$
the sequence
$v^{\lambda}(\,\cdot\,;\alpha_n,\tau_n,\vartheta_1,\vartheta_2)$
converges to $v(y)$ in $H^1_{\text{loc}}(\Omega\setminus\sigma)$.

The next point of the proof is that the right-hand side of (5.1)
coincides with $G_{\lambda}(y)$.  This is a direct computation.
The second term of the right hand side of (5.1) becomes
$$
\displaystyle
\frac{\lambda}{8\pi}\,\int_0^1J_1(\lambda\vert y-x\vert\sqrt{1-t})\,\frac{dt}{\sqrt{t(1-t)}}\\
=\frac{\lambda}{4\pi}\int_0^1(1-w^2)^{-1/2}J_1(\lambda\vert y-x\vert w)dw.
\tag {5.9}
$$
Here we know that
$$\displaystyle
s\int_0^1(1-w^2)^{-1/2}J_1(sw)dw
=1-\cos\,s.
\tag {5.10}
$$
This is an implication of (30) for $\beta=0$, $a=s$ and $\nu=1$ on p.337 of \cite{B2}:
$$\displaystyle
\int_0^s(s^2-x^2)^{-1/2}J_1(x)dx
=\frac{\pi}{2}\{J_{1/2}(\frac{s}{2})\}^2
$$
where $J_{1/2}$ stands for the Bessel function of order $1/2$
and has the explicit form
$$\displaystyle
J_{1/2}(w)=\sqrt{\frac{2}{\pi\,w}}\sin\,w,\,\,w>0.
$$
Using (5.1), (5.9) and (5.10), we conclude that
$v(y)=G_{\lambda}(y)$.

\noindent
$\Box$

{\bf\noindent Definition 5.3.}
We call the needle sequence $\{v^{\lambda}(\,\cdot\,-x;\alpha_n,\tau_n,\vartheta_1,\vartheta_2)\vert_{\Omega}\}$ given
in Theorem 5.1 a {\it standard} needle sequence for $(x,\sigma)$ for the Helmholtz equation
$\triangle u+\lambda^2 u=0$.

Therefore, now we can say that the probe method applied to the inverse
obstacle scattering problems becomes a completely {\it explicit method}
if one uses only straight needles and standard needle sequences.
In this case this explicit method gives information
on the location and shape of unknown obstacles more than the convex hull.

Some remarks are in order.

{\bf\noindent Remark 5.2.}
From Theorem 4.3 and (5.10) we get the precise values of $v^{\lambda}(y-x;\alpha,\tau,\vartheta_1,\vartheta_2)$
on the line $y=x+\,s\,\omega (-\infty<s<\infty)$ with $\omega=\vartheta_1\times\vartheta_2$.  The results are:

if $y=x+s\,\omega$ with $s\not=0$, then
$$\begin{array}{c}
\displaystyle
v^{\lambda}(y-x;\alpha,\tau,\vartheta_1,\vartheta_2)=
\frac{1}{4\pi}
\frac{E_{\alpha}(\tau s)-\cos\,\lambda s}
{s}\\
\\
\displaystyle
-\frac{\lambda}{4\pi}
\int_0^1(1-w^2)^{-1/2}E_{\alpha}(\tau(1-w^2)s)
J_1(\lambda s w)dw;
\end{array}
$$

if $y=x$, then
$$\displaystyle
v^{\lambda}(y-x;\alpha,\tau,\vartheta_1,\vartheta_2)\vert_{y=x}
=\frac{\tau}{4\pi\Gamma(1+\alpha)}.
$$

\noindent
Moreover, from (4.16) and the fact that $J_1(s)$ is an odd function we see that
$\nabla v^{\lambda}(y-x;\alpha,\tau,\vartheta_1,\vartheta_2)$ on the line
$y=x+s\,\omega\,(-\infty<s<\infty)$ is also parallel to $\omega$.  In particular,
we have
$$\displaystyle
\nabla v^{\lambda}(y-x;\alpha,\tau,\vartheta_1,\vartheta_2)\vert_{y=x}
=\frac{\tau^2}{4\pi\,\Gamma(1+2\alpha)}\,\omega.
$$

{\bf\noindent Remark 5.3.}
Since the function
$$\displaystyle
i\,\frac{\sin\,\lambda\vert y\vert}{4\pi\vert y\vert}
$$
satisfies the Helmholtz equation in the whole space, we see that the sequence
$\{\tilde{v}^{\lambda}(y-x;\alpha_n,\tau_n,\vartheta_1,\vartheta_2)\}$ given by
$$\displaystyle
v^{\lambda}(y-x;\alpha_n,\tau_n,\vartheta_1,\vartheta_2)
+i\,\frac{\sin\,\lambda\vert y-x\vert}{4\pi\vert y-x\vert},\,\,n=1,\cdots,\, y\in\Omega
$$
is also a needle sequence for $(x,\sigma)$ with $G=\tilde{G_{\lambda}}$ given by
$$\displaystyle
\tilde{G_{\lambda}}(y)=\frac{e^{i\lambda\vert y\vert}}{4\pi\vert y\vert}.
$$
This $\tilde{G_{\lambda}}$ is nothing but the fundamental solution for the Helmholtz equation
and satisfies (2.1) for any finite cone $V$ with vertex at $0$.  We also call the sequence
$\{\tilde{v}^{\lambda}(y-x;\alpha_n,\tau_n,\vartheta_1,\vartheta_2)\}$
the {\it standard needle sequence} for $(x,\sigma)$.

\begin{figure}[htbp]
\begin{center}
\epsfxsize=10cm
\epsfysize=10cm
\epsfbox{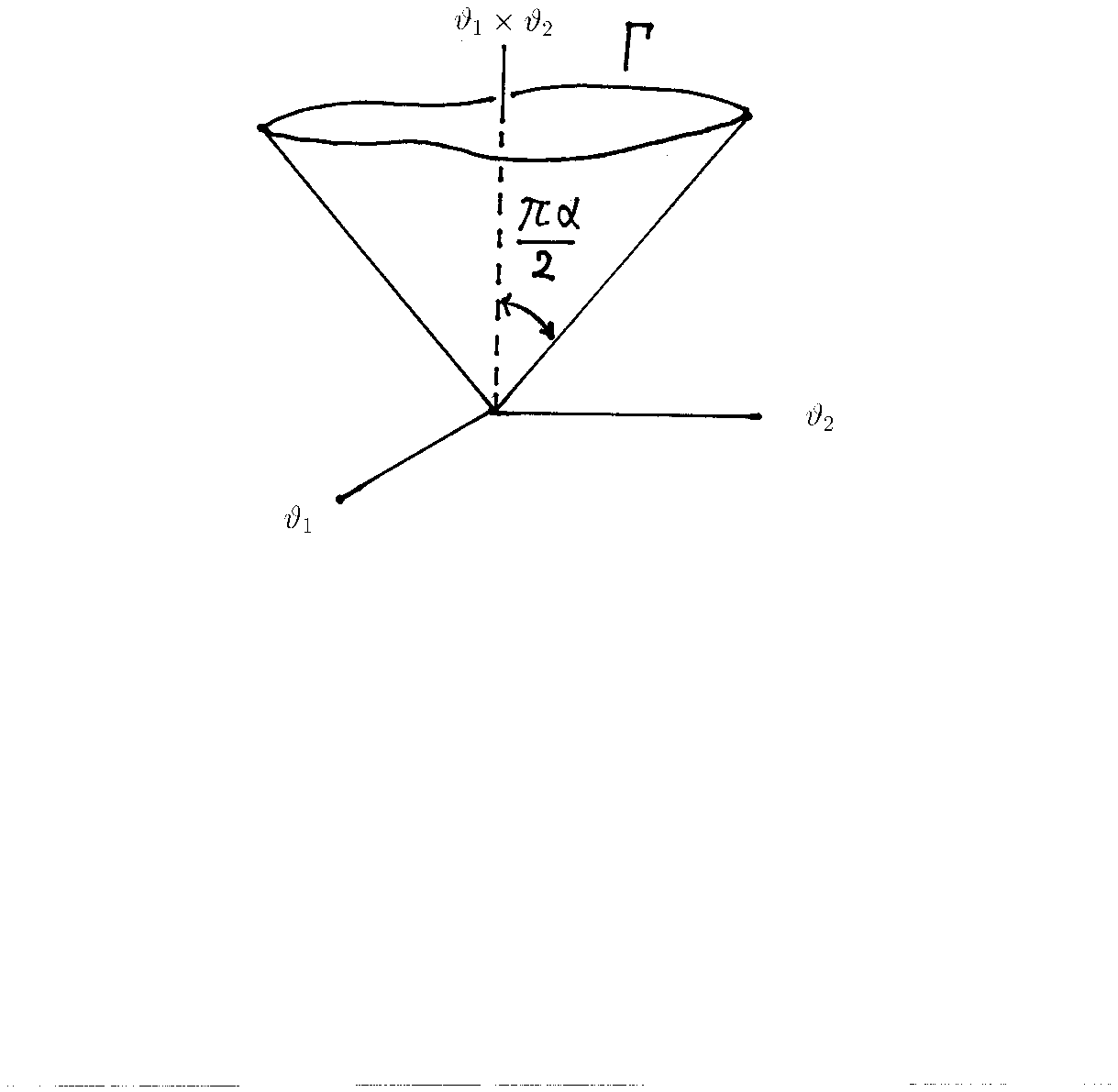}
\caption{An illustration of $\Omega$.
}\label{fig2}
\end{center}
\end{figure}

{\bf\noindent Remark 5.4.}
Fix $\alpha\in\,]0,\,1[$.
As a corollary of the proof of Theorem 5.1 we see that the function
$$
\Phi(y)=G_{\lambda}(y)-v^{\lambda}(y;\alpha,\tau,\vartheta_1,\vartheta_2),\,\,y\not=0
$$
with a large parameter $\tau$ is a Carleman function (for the
Helmholtz equation) for a special domain $\Omega$ whose boundary
consists of a part of the conic surface and a smooth surface
$\Gamma$ lying inside the cone of axis direction
$\vartheta_1\times\vartheta_2$ and aperture angle $\pi\alpha$.
See Figure \ref{fig2} for an illustration of $\Omega$. This is an
extension of a result in \cite{YM} to the Helmholtz equation and
new.

\newpage

\vspace{-1.0cm}

As pointed out in the
beginning of this section Theorem 4.1 for $\lambda>0$ does not
give this result.

{\bf\noindent Remark 5.5.}
For a recent application of the Vekua
transform to a boundary value problem for the Helmholtz equation
see \cite{CD}.  In \cite{IH}, we applied the idea of the Vekua transform
to the construction of the Herglotz wave function that approximates in a disc
a complex exponential solution of the Helmholtz equation.  This is an application
to inverse obstacle scattering problems in two dimensions.

\section{Conclusion and open problems}

We could find an explicit needle sequence for a needle given by a
segment in three dimensions.  This makes the probe method
completely explicit as same as the enclosure method. It is based
on a principle that connects the regular part of a Carleman
function with a needle sequence in the probe method.

Some further problems to be solved are in order.

\noindent $\bullet$ A mathematically interesting question is: can
one give an explicit needle sequence having a closed-form for a general needle
not necessary given by a segment in three dimensions?

\noindent $\bullet$ Using the standard needle sequence for the
Helmholtz equation $\triangle u+\lambda^2u=0$, one can construct
the corresponding indicator sequence for some inverse boundary
value problems that are the reduction of inverse obstacle
scattering problems at a fixed wave number to a bounded domain.
Can one prove the blowing up of the indicator sequence {\it
without assuming} that $\lambda$ is {\it small} when the tip of
the needle is located on the boundary, inside or passing through
unknown obstacles.  For the case when $\lambda$ is small see
\cite{Ik12,IkIm}.

\noindent $\bullet$  Construct a needle sequence for the
stationary Schr\"odinger equation $-\triangle u+V(x)u=0$
or the equation $\nabla\cdot\gamma(x)\nabla u=0$ as explicit as possible.
This is important for the
reconstruction problem of discontinuity embedded in an inhomogeneous
conductive medium.  See \cite{Ii} for the probe method applied to this problem.

\noindent
$\bullet$  Do the numerical implementation of the probe method by using the standard needle sequences
and compare the results with the previous numerical study done in \cite{CLN, EP}
of the probe method in two dimensions.

$$\quad$$

\centerline{{\bf Acknowledgement}}

This research was partially supported by Grant-in-Aid for
Scientific Research (C)(No.  18540160) of Japan Society for the
Promotion of Science.  The author wishes to thank the anonymous
referees for their comments and suggestions.

\appendix

\section{Appendix.  Proof of Theorem 4.1}

\subsection{Preliminaries}

Given an entire function $K(w)$ set
$$\displaystyle
\Phi(x)=\int_0^{\infty}\frac{K(w)}{w}\frac{\Psi(\lambda u)}{\sqrt{\vert x'\vert^2+u^2}}du
$$
where $\Psi$ is a smooth function on $[0,\,\infty[$ and
$$
\displaystyle
w=x_3+i\sqrt{\vert x'\vert^2+u^2}.
$$
Set
$$\displaystyle
f(w)=\frac{K(w)}{w}.
$$
Here we show that under suitable growth conditions for $f(w)$ as $u\longrightarrow\infty$ and
the choice of $\Psi(\lambda u)$
the function $\Phi(x)$ satisfies
$$\displaystyle
\triangle\Phi+\lambda^2\Phi=0
$$
for $x'\not=0$.

First we do the following formal computations.
We have
$$\displaystyle
\frac{\partial^2\Phi}{\partial x_3^2}
=\int_0^{\infty}f^{''}(w)\frac{\Psi(\lambda u)}{\sqrt{\vert x'\vert^2+u^2}}du.
\tag {A.1}
$$
Let $k=1,2$.  We have
$$\displaystyle
\frac{\partial\Phi}{\partial x_k}
=\int_0^{\infty}f'(w)i\frac{x_k}{\vert x'\vert^2+u^2}\Psi(\lambda u)du
-\int_0^{\infty}f(w)\frac{x_k}{(\sqrt{\vert x'\vert^2+u^2})^3}\Psi(\lambda u)du.
$$
Then we obtain
$$\begin{array}{l}
\displaystyle
\frac{\partial^2\Phi}{\partial x_k^2}
=-\int_0^{\infty}f^{''}(w)\frac{x_k^2}{(\sqrt{\vert x'\vert^2+u^2})^3}\Psi(\lambda u)du
+\int_0^{\infty}f'(w)i\frac{\Psi(\lambda u)du}{\vert x'\vert^2+u^2}\\
\\
\displaystyle
-\int_0^{\infty}f'(w)i\frac{2x_k^2}{(\vert x'\vert^2+u^2)^2}\Psi(\lambda u)du
-\int_0^{\infty}f'(w)i\frac{x_k^2}{(\vert x'\vert^2+u^2)^2}\Psi(\lambda u)du\\
\\
\displaystyle
-\int_0^{\infty}f(w)\frac{\Psi(\lambda u)du}{(\sqrt{\vert x'\vert^2+u^2})^3}
+\int_0^{\infty}f(w)3\frac{x_k^2}{(\sqrt{\vert x'\vert^2+u^2})^5}\Psi(\lambda u)du.
\end{array}
$$
This yields
$$\begin{array}{l}
\displaystyle
\sum_{k=1}^{2}
\frac{\partial^2\Phi}{\partial x_k^2}
=-\int_0^{\infty}f^{''}(w)\frac{\vert x'\vert^2}{(\sqrt{\vert x'\vert^2+u^2})^3}
\Psi(\lambda u)du
+2i\int_0^{\infty}f'(w)\frac{\Psi(\lambda u)du}{\vert x'\vert^2+u^2}\\
\\
\displaystyle
-2i\int_0^{\infty}f'(w)\frac{\vert x'\vert^2}{(\vert x'\vert^2+u^2)^2}\Psi(\lambda u)du
-i\int_0^{\infty}f'(w)\frac{\vert x'\vert^2}{(\vert x'\vert^2+u^2)^2}\Psi(\lambda u)du\\
\\
\displaystyle
-2\int_0^{\infty}f(w)\frac{\Psi(\lambda u)du}{(\sqrt{\vert x'\vert^2+u^2})^3}
+3\int_0^{\infty}f(w)\frac{\vert x'\vert^2}{(\sqrt{\vert x'\vert^2+u^2})^5}\Psi(\lambda u)du\\
\\
\displaystyle
=\int_0^{\infty}f(w)\{\frac{3\vert x'\vert^2}{(\sqrt{\vert x'\vert^2+u^2})^5}
-\frac{2}{(\sqrt{\vert x'\vert^2+u^2})^3}\}\Psi(\lambda u)du\\
\\
\displaystyle
+i\int_0^{\infty}f'(w)\{\frac{2}{\vert x'\vert^2+u^2}-\frac{3\vert x'\vert^2}{(\vert x'\vert^2+u^2)^2}\}
\Psi(\lambda u)du\\
\\
\displaystyle
-\int_0^{\infty}f^{''}(w)\frac{\vert x'\vert^2}{(\sqrt{\vert x'\vert^2+u^2})^3}\Psi(\lambda u)du\\
\\
\displaystyle
=
\int_0^{\infty}f(w)\frac{\vert x'\vert^2-2u^2}{(\sqrt{\vert x'\vert^2+u^2})^5}\Psi(\lambda u)du
+i\int_0^{\infty}f'(w)\frac{2u^2-\vert x'\vert^2}{(\vert x'\vert^2+u^2)^2}\Psi(\lambda u)du\\
\\
\displaystyle
-\int_0^{\infty}f^{''}(w)\frac{\vert x'\vert^2}{(\sqrt{\vert x'\vert^2+u^2})^3}\Psi(\lambda u)du.
\end{array}
\tag {A.2}
$$

\noindent
From (A.1) and (A.2) we have
$$\begin{array}{l}
\displaystyle
\triangle\Phi
=\int_0^{\infty}f(w)\frac{\vert x'\vert^2-2u^2}{(\sqrt{\vert x'\vert^2+u^2})^5}\Psi(\lambda u)du
+i\int_0^{\infty}f'(w)\frac{2u^2-\vert x'\vert^2}{(\vert x'\vert^2+u^2)^2}\Psi(\lambda u)du\\
\\
\displaystyle
+\int_0^{\infty} f^{''}(w)\{\frac{1}{\sqrt{\vert x'\vert^2+u^2}}
-\frac{\vert x'\vert^2}{(\sqrt{\vert x'\vert^2+u^2})^3}\}\Psi(\lambda u)du\\
\\
\displaystyle
=\int_0^{\infty}f(w)\frac{\vert x'\vert^2-2u^2}{(\sqrt{\vert x'\vert^2+u^2})^5}\Psi(\lambda u)du
+i\int_0^{\infty}f'(w)\frac{2u^2-\vert x'\vert^2}{(\vert x'\vert^2+u^2)^2}\Psi(\lambda u)du\\
\\
\displaystyle
+\int_0^{\infty}f^{''}(w)\frac{u^2}{(\sqrt{\vert x'\vert^2+u^2})^3}\Psi(\lambda u)du.
\end{array}
\tag {A.3}
$$
Integration by parts yields
$$\begin{array}{l}
\displaystyle
\int_0^{\infty}f^{''}(w)\frac{u^2}{(\sqrt{\vert x'\vert^2+u^2})^3}\Psi(\lambda u)du
=-i\int_0^{\infty}\frac{\partial}{\partial u}\{f'(w)\}
\frac{u}{\vert x'\vert^2+u^2}\Psi(\lambda u)du\\
\\
\displaystyle
=-if'(w)\frac{u}{\vert x'\vert^2+u^2}\Psi(\lambda u)\vert_0^{\infty}
+i\int_0^{\infty}f'(w)\frac{\partial}{\partial u}
\{\frac{u\Psi(\lambda u)}{\vert x'\vert^2+u^2}\}du\\
\\
\displaystyle
=i\int_0^{\infty}f'(w)\frac{\vert x'\vert^2-u^2}{(\vert x'\vert^2+u^2)^2}\Psi(\lambda u)du
+i\lambda\int_0^{\infty}f'(w)\frac{u\Psi'(\lambda u)}{\vert x'\vert^2+u^2}du.
\end{array}
$$
Note that we assumed that
$$\displaystyle
\lim_{u\longrightarrow\infty}f'(w)\frac{u}{\vert x'\vert^2+u^2}\Psi(\lambda u)=0.
$$
In what follows we always assume this type of growth condition for $f(w)$ to delete
a term coming from integration by parts.

Therefore (A.3) becomes
$$\begin{array}{l}
\displaystyle
\triangle\Phi
=\int_0^{\infty}f(w)\frac{\vert x'\vert^2-2u^2}{(\sqrt{\vert x'\vert^2+u^2})^5}\Psi(\lambda u)du\\
\\
\displaystyle
+i\int_0^{\infty}f'(w)\frac{u^2\Psi(\lambda u)}{(\vert x'\vert^2+u^2)^2}du
+i\lambda\int_0^{\infty}f'(w)\frac{u\Psi'(\lambda u)}{\vert x'\vert^2+u^2}du.
\end{array}
$$

\noindent
Integration by parts yields
$$\begin{array}{l}
\displaystyle
i\int_0^{\infty}f'(w)\frac{u^2\Psi(\lambda u)}
{(\vert x'\vert^2+u^2)^2}
=\int_0^{\infty}\frac{\partial }{\partial u}\{f(w)\}
\frac{u\Psi(\lambda u)}{(\sqrt{\vert x'\vert^2+u^2})^3}du\\
\\
\displaystyle
=f(w)\frac{u\Psi(\lambda u)}{(\sqrt{\vert x'\vert^2+u^2})^3}\vert_0^{\infty}
-\int_0^{\infty}f(w)\frac{\partial}{\partial u}\{\frac{u\Psi(\lambda u)}
{(\sqrt{\vert x'\vert^2+u^2})^3}\}du\\
\\
\displaystyle
=-\int_0^{\infty}f(w)\frac{\vert x'\vert^2-2u^2}{(\sqrt{\vert x'\vert^2+u^2})^5}
\Psi(\lambda u)du
-\lambda \int_0^{\infty}f(w)\frac{u\Psi'(\lambda u)}{(\sqrt{\vert x'\vert^2+u^2})^3}du.
\end{array}
$$
Therefore we obtain
$$\displaystyle
\triangle\Phi
=-\lambda\int_0^{\infty}f(w)\frac{u\Psi'(\lambda u)}{(\sqrt{\vert x'\vert^2+u^2})^3}du
+i\lambda\int_0^{\infty}f'(w)\frac{u\Psi'(\lambda u)}{\vert x'\vert^2+u^2}du.
$$
Here we assume that
$$\displaystyle
\Psi'(0)=0.
\tag {A.4}
$$
Then integration by parts yields
$$\begin{array}{l}
\displaystyle
i\int_0^{\infty}f'(w)\frac{u\Psi'(\lambda u)}{\vert x'\vert^2+u^2}du
=\int_0^{\infty}\frac{\partial}{\partial u}
\{f(w)\}
\frac{\Psi'(\lambda u)}{\sqrt{\vert x'\vert^2+u^2}}du\\
\\
\displaystyle
=f(w)\frac{\Psi'(\lambda u)}{\sqrt{\vert x'\vert^2+u^2}}\vert_0^{\infty}
-\int_0^{\infty}f(w)\frac{\partial}{\partial u}
\frac{\Psi'(\lambda u)}{\sqrt{\vert x'\vert^2+u^2}}du\\
\\
\displaystyle
=\int_0^{\infty}f(w)\frac{u\Psi'(\lambda u)}{(\sqrt{\vert x'\vert^2+u^2})^3}du
-\lambda\int_0^{\infty}f(w)\frac{\Psi^{''}(\lambda u)}{\sqrt{\vert x'\vert^2+u^2}}du
\end{array}
$$
and thus we obtain
$$\displaystyle
\triangle\Phi=-\lambda^2\int_0^{\infty}f(w)\frac{\Psi^{''}(\lambda u)}
{\sqrt{\vert x'\vert^2+u^2}}du.
\tag {A.5}
$$
Therefore if $\Psi$ satisfies (A.4) and
$$\displaystyle
\Psi^{''}=\Psi,
\tag {A.6}
$$
then from (A.5) we obtain
$$\displaystyle
\triangle\Phi+\lambda^2\Phi=0
$$
for $x'\not=0$.

The $\Psi$ that satisfies (A.4) and (A.6) is given by
$$
C(e^{u}+e^{-u})
$$
and thus we choose
$$
\Psi(\lambda u)=C(e^{\lambda u}+e^{-\lambda u})
$$
where $C$ is a constant.
Then it is easy to check that the formal computations given above can be justified
under the growth condition (4.1).

{\bf\noindent Remark A.1.}
We give some historical remarks.
In \cite{Y1} he considered the case when $\lambda=0$.  Therein he gave a proof
for the fact that $\Phi$ is harmonic for $x'\not=0$.
However, the proof
uses a change of variables and
can not cover the case when $K(w)$ satisfies (4.1) with $\lambda=0$,
for example, $K(w)=e^{m w}$ with $m>0$ as pointed
out in Appendix A of \cite{Iu}.  However, in \cite{Y5}, he had given a direct proof
in that case.  Using the idea of the proof, one can easily obtain that,
in the case when $\lambda=0$ $\Phi$ is harmonic for $x'\not=0$ provided
$$
\forall R>0\,\,\sup_{\vert\text{Re}\,w\vert<R}\vert K^{(m)}(w)\vert<\infty
$$
for $m=0, 1, 2$.

\noindent
Except for the explanation for the choice of a suitable $\Psi$
the proof presented in this subsection follows
that of \cite{Y5}.

\subsection{Extracting singularity at $x=0$}

We assume that $K(w)$ is real for real $w$.

\noindent
Set
$$\displaystyle
\tilde{\Phi}(x)
=\int_0^{\infty}\text{Im}\,(\frac{K(w)}{w})\frac{\Psi(\lambda u)}{\sqrt{\vert x'\vert^2+u^2}}du.
$$
Here we extract the singularity of this function at $x=0$.

Since
$$\displaystyle
K(w)=\sum_{n=0}^{\infty}\frac{K^{(n)}(0)}{n!}w^n,
$$
we have
$$\begin{array}{l}
\displaystyle
\frac{K(w)}{w}
=\frac{K(0)}{x_3+i\sqrt{\vert x'\vert^2+u^2}}+K'(0)+
+\sum_{n=2}^{\infty}\frac{K^{(n)}(0)}{n!}
(x_3+i\sqrt{\vert x'\vert^2+u^2})^{n-1}\\
\\
\displaystyle
=\frac{K(0)(x_3-i\sqrt{\vert x'\vert^2+u^2})}
{\vert x\vert^2+u^2}+K'(0)\\
\\
\displaystyle
+\sum_{n=1}^{\infty}\frac{K^{(n+1)}(0)}{(n+1)!}
\sum_{j=0}^n\frac{n!}{j!(n-j)!}x_3^{n-j}i^j(\vert x'\vert^2+u^2)^{j/2}.
\end{array}
$$
All coefficients $K^{(n)}(0), n=0, 1,\cdots$ are real and this yields
$$\begin{array}{l}
\displaystyle
\text{Im}\,(\frac{K(w)}{w})
=-\frac{\sqrt{\vert x'\vert^2+u^2}}{\vert x\vert^2+u^2}K(0)\\
\\
\displaystyle
+\sum_{n=1}^{\infty}\frac{K^{(n+1)}(0)}{(n+1)!}
\sum_{\text{$j$ is odd}}^n \frac{n!}{j!(n-j)!}
x_3^{n-j}(-1)^{(j-1)/2}(\vert x'\vert^2+u^2)^{j/2}.
\end{array}
$$
Therefore we obtain
$$\begin{array}{l}
\displaystyle
\int_0^1\text{Im}\,(\frac{K(w)}{w})\frac{\Psi(\lambda u)}{\sqrt{\vert x'\vert^2+u^2}}du
=-K(0)\int_0^1\frac{\Psi(\lambda u)}{\vert x\vert^2+u^2}du\\
\\
\displaystyle
+\sum_{n=1}^{\infty}
\frac{K^{(n+1)}(0)}{(n+1)!}
\sum_{\text{$j$ is odd}}^{n}\frac{n!}{j!(n-j)!}
x_3^{n-j}(-1)^{(j-1)/2}\int_0^1(\vert x'\vert^2+u^2)^{(j-1)/2}\Psi(\lambda u)du\\
\\
\displaystyle
\equiv
-K(0)\int_0^1\frac{\Psi(\lambda u)}{\vert x\vert^2+u^2}du
\end{array}
$$
modulo a $C^{\infty}$ function on the whole space.
Since
$$\displaystyle
\int_1^{\infty}\text{Im}\,(\frac{K(w)}{w})\frac{\Psi(\lambda u)}{\sqrt{\vert x'\vert^2+u^2}}du
$$
is $C^2$ on the whole space provided $K(w)$ satisfies a suitable growth condition,
we obtain that
$$\displaystyle
\tilde{\Phi}(x)\equiv -K(0)\int_0^1\frac{\Psi(\lambda u)}{\vert x\vert^2+u^2}du
\tag {A.7}
$$
modulo a $C^2$ function on the whole space.

Now we take
$$\displaystyle
\Psi(\lambda u)=(e^{\lambda u}+e^{-\lambda u})/2.
$$
Let $0<\vert x\vert<1$.
Then
$$\begin{array}{l}
\displaystyle
\int_0^1\frac{\Psi(\lambda u)}{\vert x\vert^2+u^2}du
=\frac{1}{2}\int_{-1}^1\frac{e^{\lambda u}}{\vert x\vert^2+u^2}du\\
\\
\displaystyle
=\pi i\text{Res}\,(\frac{e^{\lambda z}}{\vert x\vert^2+z^2})\vert_{z=i\vert x\vert}
-\frac{1}{2}\int_{z=e^{i\theta}, \,0\le\theta\le\pi}\frac{e^{\lambda z}}{\vert x\vert^2+z^2}dz\\
\\
\displaystyle
=\frac{\pi}{2}\frac{e^{i\lambda\vert x\vert}}{\vert x\vert}
-\frac{i}{2}\int_0^{\pi}\frac{e^{\lambda(\cos\,\theta+i\sin\,\theta)}}
{\vert x\vert^2+e^{2i\theta}}e^{i\theta}d\theta
\equiv
\frac{\pi}{2}\frac{e^{i\lambda\vert x\vert}}{\vert x\vert}
\end{array}
$$
modulo a $C^{\infty}$ function on $\vert x\vert<1$.  Therefore from (A.7)
one concludes that
$$\displaystyle
\tilde{\Phi}(x)
\equiv -\frac{K(0)\pi}{2}\frac{e^{i\lambda\vert x\vert}}{\vert x\vert}
$$
modulo a $C^2$ function on the whole space provided $K(w)$ satisfies (4.1).
In particular, if $K(0)=1$, then
$$\displaystyle
-\frac{1}{2\pi^2}\tilde{\Phi}(x)
\equiv\frac{e^{i\lambda\vert x\vert}}{4\pi\vert x\vert}
$$
modulo a $C^2$ function on the whole domain.

Using the results in subsections, one obtains Theorem
4.1.

\vskip1cm
\noindent
e-mail address

ikehata@math.sci.gunma-u.ac.jp
\end{document}